\documentclass[10pt,reqno]{article}
\usepackage{amsmath, amsthm, amssymb, stmaryrd}
\usepackage{hyperref}
\usepackage{enumerate}
\usepackage{url}
\usepackage{color}
\usepackage{tikz}
\usetikzlibrary{patterns}
\usepackage[T1,T2A]{fontenc}
\usepackage[utf8]{inputenc}

\usepackage{amsthm}

\newtheorem{thm}{Theorem}
\newtheorem{lemma}{Lemma}

\theoremstyle{definition}

\theoremstyle{remark}

\def \le {\leqslant}
\def \ge {\geqslant}

\usepackage{tikz} 
\usetikzlibrary{calc,arrows,decorations.pathreplacing,fadings,3d,positioning}

\usepackage[mathscr]{euscript}

\begin{document}

\begin{Huge}
\centerline{Approximation to a real number and its cube }
\centerline{ by rational numbers}
\end{Huge}
\vskip+0.5cm
\begin{large}
\centerline{Nikolay Moshchevitin}
\end{large}

\vskip+1cm

\section{Introduction and main results}

\subsection{Irrationality measure functions and Diophantine exponents}

For a real vector $\pmb{\xi} = (\xi_1,\xi_2) $ we consider the problem of simultaneous approximation by rationals and the dual problem dealing with linear form. Let
$||\cdot ||$ be the distance to the nearest integer and 
\begin{equation}\label{irrat}
   \psi_{\pmb{\xi}}(t) = \min_{\begin{array}{c}q\in \mathbb{Z}
 :\cr 0<q \le t \end{array}}
  \max_{j=1,2} ||q\xi_j||,\,\,\,\,\,\,\,\,\,
   \psi^*_{\pmb{\xi}}(t) = \min_{\begin{array}{c}(m_1,m_2) \in \mathbb{Z}^2:\cr 0<\max(|m_1|,|m_2|) \le t \end{array}}
   ||m_1\xi_1+m_2\xi_2||
   \end{equation}
   be irrationality measure functions related to simultaneous approximation and linear form setting correspondingly.
   Vector
   $\pmb{\xi}= (\xi_1,\xi_2) $  is called {\it singular} if 
   $$
   \lim_{t\to \infty} t^{1/2} \psi_{\pmb{\xi}}(t) = 0
   \,\,\,
   \,\,\,\,\,\text{
   or equivalently}\,\,\,\,\,\,\,\,
    \lim_{t\to \infty} t^{2} \psi_{\pmb{\xi}}^*(t) = 0
   $$
   (see \cite{khi}).
   Uniform Diophantine exponents  for these two settings are defined as
\begin{equation}\label{expo}
   \hat{\omega} ({\pmb{\xi}}) =
   \sup\left\{\gamma\in \mathbb{R}:\,\,
   \limsup_{t\to \infty} t^\gamma \cdot \psi_{\pmb{\xi}} (t)<\infty\right\},\,\,\,\,\,\,\,\,
   \text{and}
   \,\,\,\,\,\,\,\,
    \hat{\omega}^* ({\pmb{\xi}}) =
   \sup\left\{\gamma\in \mathbb{R}:\,\,
   \limsup_{t\to \infty} t^\gamma \cdot \psi_{\pmb{\xi}}^* (t)<\infty 
   \right\}.
   \end{equation}
     It is well known that for $\xi_1,\xi_2$ linearly independent over $ \mathbb{Q}$ together with 1 Diophantine exponents satisfy inequalities 
  $$
  \frac{1}{2} \le    \hat{\omega} ({\pmb{\xi}}) \le 1,\,\,\,\,\,\,\, 2\le    \hat{\omega}^* ({\pmb{\xi}})б
  $$ 
  and the famous Jarn\'{\i}k's equality \cite{j} states that 
   \begin{equation}\label{ja}
   \hat{\omega}({\pmb{\xi}}) + \frac{1}{\hat{\omega}^*({\pmb{\xi}}) } = 1.
      \end{equation}

\subsection{Uniform approximation to real number and its cube}

In \cite{Roy}  the authors considered the case
$$
\xi_1 = \xi,\,\,\,\,\, \xi_2 = \xi^3.
$$
In this case for the irrationality functions (\ref{irrat}) we will simply use notation
$$
\psi(t) = \psi_{(\xi,\xi^3)}(t),\,\,\,\,\,\,\,\,\,
\psi^*(t) = \psi^*_{(\xi,\xi^3)}(t)
$$
and for exponents (\ref{expo}) also 
$$
\hat{\omega} = \hat{\omega}((\xi,\xi^3)),\,\,\,\,\,\,\,\,
\hat{\omega}^* = \hat{\omega}^*((\xi,\xi^3)),
$$
so in particular
$$
   \psi^*(t) = \min_{\begin{array}{c}(m_1,m_2) \in \mathbb{Z}^2:\cr 0<\max(|m_1|,|m_2|) \le t \end{array}}
   ||m_1\xi+m_2\xi^3||\,\,\,\,\,\,\, \text{and}\,\,\,\,\,\,\,
   \hat{\omega}^*  =
   \sup\left\{\gamma\in \mathbb{R}:\,\,
   \limsup_{t\to \infty} t^\gamma \cdot \psi^* (t)<\infty
  \right\}.
   $$

 In this notation the main result of \cite{Roy} can be stated as follows:
\begin{equation}\label{equa0}
 \text{\it inequality }\,\,\,  \hat{\omega} \le 
 \frac{2(9+\sqrt{11})}{35} = 0.7038^+
 \,\,\,\text{\it holds in the case when} \,\,\,1, \xi,\xi^3\,\,\,
 \text{\it are linearly independent over }\,\,\,\mathbb{Q}.
\end{equation}

In the present paper we  prove the following 

\begin{thm}\label{T0}
There exists transcendental $\xi$ such that $ \hat{\omega}^* \ge
\frac{4+\sqrt{7}}{3} =  2.2152^+.
 $
\end{thm}

 Jarn\'{\i}k's equality gives alternative formulation of Theorem  \ref{T0}:
\begin{equation}\label{equa}
\text{\it 
There exists transcendental}\,\,\,\xi \,\,\,\text{such that}\,\,\, \hat{\omega} \ge \frac{\sqrt{7}-1}{3} = 0.5485^+.
\end{equation} 
The question if there exists transcendental $\xi$ such that  $ \hat{\omega} \in \left(  \frac{\sqrt{7}-1}{3} , \frac{2(9+\sqrt{11})}{35} \right]$ remains open. In \cite{Roy} the authors claim that "{\it the method of their paper 
is capable of lowering the upper bound in (\ref{equa0}), possibly down to 
$\frac{1+3\sqrt{5}}{11} = 0.7007^+$}",
 and we completely agree with this claim.

\subsection{General result  about polynomial curve}

Let $d\ge 3$ be an integer. For a real $\tau>2d-2$ we consider four functions
\begin{equation}\label{g}
g_1(\tau)=\frac{\tau^2}{(d-1)(\tau+1)},\,\,\,\,\,\,g_2(\tau)=\frac{\tau+1}{d},\,\,\,\,\,\,g_3(\tau)=2+\frac{d-2}{\tau-d+2},\,\,\,\,\,\,g_4(\tau)=2+\frac{2}{\tau}-\frac{d-1}{\tau^2},
\end{equation}
and put
\begin{equation}\label{om}
\omega_d(\tau)=\min\left\{g_1(\tau),g_2(\tau),g_3(\tau),g_4(\tau)\right\},
\,\,\,\,\,\,\,\,\,\,\,\,
\omega_d=\sup_{\tau>2d-2}\omega_d(\tau).
\end{equation}
We distinguish cases $d = 3$ and $d\ge 4$.  In these cases calculations lead to the following conclusion.

\vskip+0.2cm
\noindent
({\bf a})
Let $d=3$. Then the supremum in {\rm(\ref{om})} is attained at the unique point $\tau=3+\sqrt7$, which is the
larger root of  polynomial
$$ Q_3 (\tau ) = \tau^2-6\tau+2,
$$ and
$
\omega_3=\frac{4+\sqrt7}{3}
$\ is the larger root of the polynomial
$$
R_3(\omega) = 3\omega^2-8\omega+3.
$$

\vskip+0.2cm
\noindent
({\bf b})
Let $d\ge 4$.
We consider two polynomials
$$
Q_d(\tau)=\tau^3-(2d-1)\tau^2-2d\tau+d(d-1),
$$
$$
R_d(\omega)=d^2\omega^3-2d(d+1)\omega^2+(2d+1)\omega+(d-1).
$$
Notice that 
$
Q_d(d\omega-1)=d\cdot R_d(\omega).
$

 The supremum in {\rm(\ref{om})} is attained at the unique point $\tau=\tau_d$, where $\tau_d$ is the
unique root of $Q_d(\tau)$ in the interval $(2d-1,2d)$; it is the largest root of $Q_d$. One has
$$
\omega_d=\frac{\tau_d+1}{d},
$$
and $\omega_d$ is the unique root of the polynomial $R_d(\omega)$ in the interval $\left(2,2+\frac1d\right)$; it is the largest
root of $R_d$.

\vskip+0.2cm

We should notice that  for $d\ge 4$ one has
$$
2+\frac{3}{4d}<\omega_d<2+\frac{3}{4d}+\frac{1}{2d^2},
$$
and the sequence 
  $\omega_3,\omega_4,\omega_5,...$ decreases to $2$.

\vskip+0.2cm
Theorem \ref{T0} is a particular case of the following general result.

 \begin{thm}\label{T2}
 Let $f(t) $ be a polynomial with rational coefficients of degree $d\ge 3$ and 
 for given $\xi$ we put 
$$
 \hat{\omega}^*[f] = \hat{\omega}^*((\xi, f(\xi))).
 $$
Then there exists transcendental $\xi$ such that 
$$
 \hat{\omega}^*[f] \ge \omega_d
 .
 $$
\end{thm}

\noindent
{\bf Remark 1.}
Our construction can be easily modified to construct an uncountable set of numbers $\xi$ satisfying the conclusion of Theorem \ref{T2}.

\noindent
 {\bf Remark 2.}
To the best of our knowledge, up to now the existence of totally irrational singular vectors
on a curve which is not a line was known only for conics defined over $\mathbb{Q}$.
For the parabola $(t,t^2)$ this is the famous result of Roy \cite{r}:
for his extremal numbers $\xi$ one has
$\hat{\omega}((\xi,\xi^2)) = \frac{\sqrt{5}-1}{2}>\frac{1}{2}$;
in \cite{rc} he extended it to all conics defined over $\mathbb{Q}$.
The general existence results from \cite{k1,k2} deal with analytic submanifolds
of dimension at least two and do not cover curves.
Theorem \ref{T2} gives totally irrational singular vectors on the curves $(t,f(t))$
for all polynomials $f$ with rational coefficients of degree $d\ge 3$.

\subsection{ Idea of the construction.}
Let $\xi$ be a number with convergents $\frac{p_\nu}{q_\nu}$ and rapidly growing denominators.
Every convergent gives a rational point
$\left(\frac{p_\nu}{q_\nu}, f\left(\frac{p_\nu}{q_\nu}\right)\right)$ of the curve $(t,f(t))$,
that is an integer vector
$$
Z_\nu = \left(q_\nu^d,\, p_\nu q_\nu^{d-1},\, q_\nu^d f\left(\frac{p_\nu}{q_\nu}\right)\right),
\,\,\,\,\,\,\,
\max\left( |q_\nu^d \xi - p_\nu q_\nu^{d-1}|,\,\left|q_\nu^d f(\xi) -q_\nu^d f\left(\frac{p_\nu}{q_\nu}\right)\right|\right) \ll \frac{q_\nu^{d-1}}{q_{\nu+1}},
$$
which is a very good simultaneous approximation to $(\xi, f(\xi))$.
In the dual setting this means the following. Consider the two-dimensional rational subspace
$$
\pi_\nu = \left\{ X = (x_0,x_1,x_2) \in \mathbb{R}^3:\,\,\, x_0+x_1\frac{p_\nu}{q_\nu}+x_2f\left(\frac{p_\nu}{q_\nu}\right) = 0\right\}
$$
orthogonal to $Z_\nu$. It contains many integer vectors $X$ with small values of
$|x_0+x_1\xi+x_2f(\xi)|$. Among them are the vectors $A_\nu$, $F_\nu$, $B_\nu$ defined below
(for a simpler proof of Section \ref{ssim} vectors $F\nu$ are not necessary).
These vectors exist for every continued fraction: no condition on partial quotients is necessary for them.
So for the linear form we automatically have a long sequence of good approximations which belong to
the same subspace $\pi_\nu$, and in $\pi_\nu$ the problem is in fact one-dimensional.

However this is not enough. Subspace $\pi_\nu$ contains no integer vector essentially better than $B_\nu$,
and every non-zero integer vector of the next subspace $\pi_{\nu+1}$ has norm at least $q_{\nu+1}$.
Between $|B_\nu|\asymp q_\nu^d$ and $|A_{\nu+1}| = q_{\nu+1}$ we need one more integer vector which
belongs neither to $\pi_\nu$ nor to $\pi_{\nu+1}$ and links these two subspaces. This is the vector $W_\nu$.
For a general continued fraction 
and general $f$ we cannot expect that
such a vector exists. The conditions on the partial quotient
$a_{\nu+1}$ in our construction have the only purpose to ensure that $W_\nu$ is an integer vector
and that the linear form is small on it. So the collections of dependent approximations
come for free from the convergents, and the special choice of partial quotients is used only to
link the successive collections one to another.

\subsection{Structure of the paper}

The paper contains two proofs based on the same idea.
Theorem \ref{T0} is the case $f(t)=t^3$ of Theorem \ref{T2}, and so formally it is enough
to prove Theorem \ref{T4}, which is a detailed version of Theorem \ref{T2}.
However the proof of the general result involves many parameters and rather cumbersome
conditions on partial quotients, and the main idea may be hidden behind the technical details.
For this reason we begin with the simplest case.
In Section 2 we consider the curve $(t,t^3)$ and prove only that the vector $(\xi,\xi^3)$ is singular
(Theorem \ref{T1}), without any attempt to get a good exponent.
This section is independent of Section 3, it is completely elementary, and we hope that
it explains the method and makes the paper easier to read.
The reader who is interested in the idea only may stop after Section 2.

Section 2 is organised as follows.
In Section 2.1 we define the number $\xi$ by conditions on the partial quotients of its continued
fraction and formulate Theorem \ref{T1}.
In Section 2.2 we show that these conditions can be satisfied.
In Section 2.3 we introduce the vectors $A_\nu$ and $B_\nu$, which are defined for every continued fraction,
and in Section 2.4 we deal with the linking vector $W_\nu$; it is the only place
where the conditions on partial quotients are used.
In Section 2.5 we combine the bounds for $A_\nu, B_\nu, W_\nu$ and finish the proof of Theorem \ref{T1}.

In Section 3 we prove Theorem \ref{T4}. The scheme of the proof is the same,
but there are several differences.

\noindent
1. For $f(t) = t^3$ the conditions on $a_{\nu+1}$ are explicit congruences.
For a general polynomial they are written in terms of the integers $\Phi_{j,\nu}$,
which play the role of derivatives of $f$ at the point $\frac{p_\nu}{q_\nu}$,
and we need a special statement (Lemma 7) to show that admissible partial quotients exist.
All these objects are introduced in Section 3.1, where Theorem \ref{T4} is formulated.

\noindent
2. For $f(t)=t^3$ the number $D_\nu$, which is responsible for the value of the linear form
at the vector $W_\nu$, can be written explicitly.
In the general case we have only an upper bound for it. It is obtained in Section 3.2.

\noindent
3. In Section 2 we need only the inequality $\psi^*(t) = o(t^{-2})$, and three vectors
$A_\nu, B_\nu, W_\nu$ for each $\nu$ are enough.
To get the exponent from Theorem \ref{T4} we need one more vector $F_\nu$,
which lies in the same two-dimensional subspace as $A_\nu$ and $B_\nu$
and improves the bound for $\psi^*(t)$ between $|A_\nu|$ and $|B_\nu|$.
Vectors $A_\nu, B_\nu, F_\nu$ are considered in Section 3.3, and vector $W_\nu$ in Section 3.4.

\noindent
4. In Section 2 the partial quotients are simply supposed to grow fast enough.
In Section 3 the growth is prescribed precisely by the parameter $\tau$,
that is $q_{\nu+1}\asymp q_\nu^\tau$, and all the bounds are given with explicit constants.
In Section 3.5 we collect the bounds for the four vectors, finish the proof of Theorem \ref{T4}
and so obtain the exponent $\omega_d(\tau)$ which should be optimised in $\tau$.

 \section{Simple construction of singular vector}\label{ssim}

\subsection{Construction of number $\xi$}

For our purpose it would be convenient to use  continued fraction
\begin{equation}\label{o}
\xi\ = \cfrac{1}{a_1+\cfrac{1}{a_2+\cfrac{1}{a_3+\cdots}}}
=[0;a_1,a_2,a_3,...]
,\,\,\,\,\, a_j \ge 1
\end{equation}
with convergents 
$\frac{p_\nu}{q_\nu}=[0;a_1,a_2,a_3,...,a_\nu]$
satisfying
\begin{equation}\label{zero}
\begin{aligned}
(p_0,q_0)&=(0,1), & \,\,\, p_{\nu+1}=a_{\nu+1}\,p_\nu+p_{\nu-1} ,&\,\,\,\,\,\,
q_\nu\,p_{\nu-1}-p_\nu\,q_{\nu-1}=(-1)^\nu,&\\
(p_1,q_1)&=(1,a_1), & \,\,\,q_{\nu+1}=a_{\nu+1}\,q_\nu+q_{\nu-1}, & \,\,\,\,\,\,\, (p_\nu,q_\nu) = 1, \,\,\,\, \nu\ge 1.&
\,\,\,\, 
\end{aligned}
\end{equation}
It is well known that 
\begin{equation}\label{11}
|q_\nu \xi- p_\nu |<\frac{1}{q_{\nu+1}}.
\end{equation}

We choose partial quotients $a_\nu$ in accordance with the following  inductive rule.  Take $a_1 = 2$.
If $a_1,...,a_\nu$ are fixed 
take  $a_{\nu+1}$ satisfying inequality
\begin{equation}\label{roth}
 a_{\nu+1}>q_\nu^{6},
\end{equation}
 which leads to 
 \begin{equation}\label{1}
 q_{\nu+1}>q_\nu^{7},
 \end{equation}
and the congruence conditions (we use Chinese Remainder Theorem)
 \begin{equation}\label{21}
a_{\nu+1}\equiv 2(-1)^{\nu+1} p_{\nu-1}q_{\nu-1}\pmod{p_\nu},
\end{equation}
\begin{equation}\label{22}
a_{\nu+1}\equiv(3p_{\nu-1}q_\nu+(-1)^\nu)\,p_{\nu-1}q_{\nu-1}\pmod{q_\nu^{2}}.
\end{equation}
Notice that in our construction $a_1 = 2 $ and $a_2 \ge 4$, and so
\begin{equation}\label{pq}
\frac{4q_\nu}{9}\leqslant p_\nu\leqslant\frac{q_\nu}{2}\,\,\,\,\,\text{for all}\,\,\,\,\, \nu\ge 1.
\end{equation}

\begin{thm}\label{T1}
Let $\xi$  be defined by continued fraction {\rm (\ref{o})} and its partial quotients satisfy 
$a_1=2$ and 
 {\rm (\ref{roth}), (\ref{21}), (\ref{22})}.
Then $\xi$ is transcendental and the  vector $(\xi,\xi^3)\in \mathbb{R}^2$ is singular.
\end{thm}

Here we should notice that non-algebraicity of $\xi$ immediately follows from Roth's theorem, because by (\ref{11}) and(\ref{roth})
inequality $\left|\xi-\frac{p}{q}\right| < \frac{1}{q^8}$ has infinitely many rational solutions $\frac{p}{q}$.

\subsection{About conditions (\ref{21}) and (\ref{22})}

Here we give some comments concerning our choice of partial quotients.

Notice that  \eqref{21} is equivalent to
 \begin{equation}\label{211}
q_{\nu+1}+q_{\nu-1}\equiv 0\pmod{p_\nu},
\end{equation}
 while \eqref{22} is equivalent to
 \begin{equation}\label{222}
p_{\nu+1}+2\,(-1)^\nu\,p_{\nu-1}^{2}\,q_\nu \equiv 0\pmod{q_\nu^{2}}.
\end{equation}
In particular, it follows from (\ref{222}) that 
 \begin{equation}\label{2222}
p_{\nu+1}\equiv 0\pmod{q_\nu}.
\end{equation}

  \begin{lemma}\label{L0}
 Three numbers  defined by equalities  
 \begin{equation}\label{wcoor}
 w_{0,\nu}=
  \frac{
 p_\nu (p_{\nu+1} + 2(-1)^\nu p_{\nu-1}^2 q_\nu)}{ q_\nu^2},\,\,\,\,
 w_{1,\nu}=
 -\frac{3(p_{\nu+1} + (-1)^\nu p_{\nu-1}^2q_\nu)}{  q_\nu},\,\,\,\,
 w_{2,\nu} =
 \frac{2(q_{\nu+1} + q_{\nu-1})}{ p_\nu} + (-1)^\nu q_{\nu-1}^2
 \end{equation}
are integers.
  \end{lemma}
  \begin{proof}
  $ w_{0,\nu}\in \mathbb{Z}$ follows from (\ref{222});
  $w_{1,\nu}\in \mathbb{Z}$ follows from (\ref{2222});
    $w_{2,\nu}\in \mathbb{Z}$ follows from (\ref{211}).
 \end{proof}

\subsection{Three-dimensional vectors $A_\nu, B_\nu$}\label{s2}
For vector $X= (x_0,x_1, x_2) \in \mathbb{R}^3$ we denote
$$
|X| = \max_{j=0,1,2}|x_j|,\,\,\,\,\, L(X) = |x_0 +x_1\xi+x_2 \xi^3|.
$$

In this section we define  integer vectors $A_\nu, B_\nu \in \mathbb{Z}^3$. Put
\begin{equation}\label{AB}
A_\nu = (-p_\nu, q_\nu, 0),
\,\,\,\,\,\,
B_\nu = (2p_\nu^{3},\ \ -3q_\nu p_\nu ^{2},\ \ q_\nu^{3}).
\end{equation}
We prove simple inequalities related to approximations with vectors  $A_\nu, B_\nu$.
 
 \begin{lemma}\label{L1}
 For any $\nu$ we have
 $$
 |A_\nu| = q_\nu,\,\,\,\,  L(A_\nu) <\frac{1}{q_{\nu+1}}.
 $$
 \end{lemma}
 \begin{proof}
  Follows immediately from (\ref{11}).
  \end{proof}
  
   \begin{lemma}\label{L2}
 For any $\nu$ we have
 $$
 |B_\nu| = q_\nu^3,\,\,\,\,  L(B_\nu) <\frac{3q_\nu}{q_{\nu+1}^2}.
 $$
 \end{lemma}
  \begin{proof}
  Equality for the norm of vector $B_\nu$ is clear. For the corresponding approximation we have
  $$
  2p_\nu^{3} -3q_\nu p_\nu ^{2}\xi + q_\nu^{3}\xi^3 = (p_\nu-q_\nu\xi)^2(q_\nu\xi +2p_\nu).
  $$
  So by   (\ref{11}) we get
  $$
  L(B_\nu)  = | 2p_\nu^{3} -3q_\nu p_\nu ^{2}\xi + q_\nu^{3}\xi^3|<
  3q_\nu |p_\nu - q_\nu \xi|^2 \le \frac{3q_\nu}{q_{\nu+1}^2}.
  $$
Lemma is proven.
  \end{proof}
 
 \subsection{Vector $W_\nu$}
 
 We need to consider one more special vector
 $$W_\nu  
 =
 xA_\nu+ y B_\nu + zA_{\nu+1}
 $$
 defined as a linear combination of vectors $ A_\nu, B_\nu , A_{\nu+1}$ with coefficients
 \begin{equation}\label{coe}
 x = (-1)^{\nu+1}\frac{6}{q_\nu^3},\,\,\,\,\,
 y = \frac{2q_{\nu+1} +q_{\nu-1}((-1)^\nu p_{\nu-1}q_\nu +1)}{p_\nu q_\nu^3},\,\,\,\,
 z = \frac{3p_\nu}{q_\nu^2}.
 \end{equation}
 We explain   properties of vector $W_\nu$ which have determined the 
  choice of coefficients (\ref{coe}).
 
 \begin{lemma}\label{L5}
 Vector $W_\nu$  has coordinates 
 $W_\nu
 =
 (w_{0,\nu}, w_{1,\nu}, w_{2,\nu}) 
 $
 defined in (\ref{wcoor}).    In particular
 $W_\nu \in \mathbb{Z}^3.$
  \end{lemma}
  \begin{proof}
  Formulas (\ref{wcoor}) for the coordinates of $W_\nu$ follow immediately from the definitions
  (\ref{AB})
   of
  $A_\nu, B_\nu, A_{\nu+1}$, formulas (\ref{zero}) and equalities for coefficients (\ref{coe}),
 and from Lemma \ref{L0} we have  $W_\nu \in \mathbb{Z}^3$.
\end{proof}

   \begin{lemma}\label{L3} For all $\nu$ large enough one has
 $q_\nu^3 \le |W_\nu| = w_{2,\nu} \le \frac{5 q_{\nu+1}}{q_\nu}$.
  \end{lemma}
  \begin{proof}
  From the growth condition (\ref{1}), inequality (\ref{pq})  and formulas (\ref{wcoor}) we see that 
  $$
  |w_{0,\nu}| \le \frac{p_{\nu+1}}{q_\nu}(1+o(1)),\,\,\,
    |w_{1,\nu}| = \frac{3p_{\nu+1}}{q_\nu} (1+o(1))\le
      \frac{3q_{\nu+1}}{2q_\nu} (1+o(1))
    ,\,\,\,\, \nu \to \infty,
  $$
  while by  (\ref{pq}) we get
  $$
      |w_{2,\nu}|  =\frac{2q_{\nu+1}}{p_\nu} (1+o(1)),\,\,\,\,
      \frac{4q_{\nu+1}}{q_\nu}(1+o(1))\leqslant w_{2,\nu}\leqslant \frac{9q_{\nu+1}}{2q_\nu}(1+o(1)),
  $$ and 
  $   |w_{2,\nu}|$
 is the maximal one.  As for upper and lower bounds for $   |w_{2,\nu}| $, they immediately follow from the 
 last formula from (\ref{wcoor}) and (\ref{1}).
\end{proof}

   \begin{lemma}\label{L4} For all $\nu$ large enough one has
 $ L(W_\nu) \le  \frac{2q_{\nu-1}^2}{q_\nu^2q_{\nu+1}^2}$.
  \end{lemma}
  \begin{proof}
  Consider  polynomials 
  $$
  A(T) = q_\nu T - p_\nu, \,\,\,
  B(T) = 2p_\nu^{3}-3q_\nu p_\nu ^{2} T + q_\nu^{3}T^3,
  \,\,\,\text{and}\,\,\,
  W(T) = w_{0,\nu} + w_{1,\nu}T +w_{2,\nu}T^3.
  $$
  Notice that 
  $$
  A\left(\frac{p_{\nu+1}}{q_{\nu+1}}\right) =\frac{(-1)^\nu }{q_{\nu+1}}\,\,\,\,
  \text{and}
  \,\,\,
  B\left(\frac{p_{\nu+1}}{q_{\nu+1}}\right) =
  \frac{ 2p_\nu^{3}q_{\nu+1}^3-3q_\nu p_\nu ^{2} p_{\nu+1}q_{\nu+1}^2 + p_{\nu+1}^3q_\nu^{3}}{q_{\nu+1}^3}=
  \frac{3p_\nu}{q_{\nu+1}^2} + \frac{(-1)^\nu}{q_{\nu+1}^3}.
  $$
  Observe that 
    $$
  W\left(\frac{p_{\nu+1}}{q_{\nu+1}}\right) = xA\left(\frac{p_{\nu+1}}{q_{\nu+1}}\right)
  + y B\left(\frac{p_{\nu+1}}{q_{\nu+1}}\right)
  ,
  $$
  and by formulas (\ref{coe}) for $x$ and $y$ after cancellation we have  
  $$
  W\left(\frac{p_{\nu+1}}{q_{\nu+1}}\right)
= (-1)^\nu\,\frac{(3p_{\nu-1}^2 q_\nu^2-1)\,q_{\nu+1}
+ q_{\nu-1}\left((-1)^\nu p_{\nu-1}q_\nu+1\right)}{p_\nu\, q_\nu^3\, q_{\nu+1}^3},
$$
and so for all $\nu$ large enough we get
    \begin{equation}\label{P2}
 \left| W\left(\frac{p_{\nu+1}}{q_{\nu+1}}\right)\right| \le \frac{7q_{\nu-1}^2}{4q_\nu^2q_{\nu+1}^2}.
  \end{equation}
  On the other hand, for all $\nu$ large enough we have
  \begin{equation}\label{P1}
 \left| W(\xi) - W\left(\frac{p_{\nu+1}}{q_{\nu+1}}\right)\right|\le
 3 |W_\nu| \cdot
  \left|\xi - \frac{p_{\nu+1}}{q_{\nu+1}}\right|\le 
  \frac{15}{q_\nu q_{\nu+2}}
 \le
  \frac{q_{\nu-1}^2}{4q_\nu^2q_{\nu+1}^2}
 .
  \end{equation}
  (we used (\ref{11}) and upper bound from Lemma \ref{L3}).
  Now combining (\ref{P2}) and (\ref{P1}) we  get inequality of Lemma \ref{L4}.
  \end{proof}
  
   \subsection{Proof of Theorem \ref{T1}}\label{sQ}
   
  For  $\nu$ large enough
 and for 
   every real $t$  from the interval  $q_\nu\le t<q_{\nu+1}$  we prove the inequality
   \begin{equation}\label{V1}
t^{2}\psi^* (t)\;<\;\frac{75}{q_\nu},
\end{equation}
which leads to
   $$
   \lim_{t\to \infty} t^2\psi^* (t) =0 ,
   $$
   and this is just singularity condition.
   
   Observe that  by equalities from Lemmas \ref{L1} and \ref{L2}  and  Lemma \ref{L3} we have
   $$ 
   q_\nu = |A_\nu|<q_\nu^3 = |B_\nu| \le |W_\nu|< q_{\nu+1} = |A_{\nu+1}|.
   $$
   For $t \in [|A_\nu|, |B_\nu| ]$ from definition of function $\psi^* (t)$,  Lemma \ref{L1}  and the  growth condition (\ref{1}) we get
   $$
   t^{2}\psi^* (t)\le
   |B_\nu| ^2L(A_\nu)\le  \frac{1}{q_{\nu}}.
   $$
      For $t \in [|B_\nu| , |W_\nu|]$ from  the inequalities   of  Lemmas \ref{L2} and \ref{L3}   we get
   $$
   t^{2}\psi^* (t)\le 
   |W_\nu| ^2L(B_\nu)
   \le  \frac{75}{q_{\nu}}.
   $$
        For $t \in [ |W_\nu|, |A_{\nu+1}|]$ from  Lemma \ref{L4}   and growth condition (\ref{1}) we get
   $$
   t^{2}\psi^* (t)\le  |A_{\nu+1}|^2L(W_\nu)\le 
   \frac{2q_{\nu-1}^2 }{q_\nu^2}\le  \frac{1}{q_{\nu}}
   $$
   for $\nu $ large enough.
   In such a way we establish (\ref{V1}).$\Box$
   \vskip+0.3cm
   
   {\bf Remark.}  By taking  in the described construction  $a_{\nu+1}\asymp q_\nu^6$ the argument of the proof of Theorem \ref{T1}  immediately gives $\hat{\omega}^*  \ge \frac{13}{6}>2$. This bound may be improved by optimising the argument of Section \ref{sQ}.
   
   \section{Proof of Theorem \ref{T2}}
   
   In our proof number $\xi$ will be also constructed as a continued fraction (\ref{o}) whose convergent fractions 
    $\frac{p_\nu}{q_\nu}$
 satisfy certain conditions related to  polynomial $f$. These conditions are rather cumbersome.
   In the next subsection we introduce all the objects which are necessary to write these conditions in terms of partial quotients $a_\nu$ of $\xi$.

   \subsection{Convergents, generalised derivatives and inductive definition of partial quotients}
   Fix $\tau >2d-2$.
   
   Without loss of generality we may assume that 
   \begin{equation}\label{fv}
f(t)=c_dt^d+c_{d-1}t^{d-1}+\dots+c_1t+c_0,\,\,\,\,\,\,\,c_j\in\mathbb{Z},\,\,\,\,\,\,\,c_d=1,\,\,\,\,\,\,\,d\ge 3.
\end{equation}
By Taylor  formula applied  for point $\frac{p_\nu}{q_\nu}$  polynomial  $f(t)$ may be written as
$$
f(t) = \sum_{j=0}^d\frac{1}{j!}f^{(j)}\left(\frac{p_\nu}{q_\nu}\right)\left(t-\frac{p_\nu}{q_\nu}\right)^j,
$$
and so
   \begin{equation}\label{fvv}
   q_\nu^df(t)=\sum_{j=0}^{d}\Phi_{j,\nu}\,(q_\nu t-p_\nu)^j ,\,\,\,\,\,
   \text{where}
   \,\,\,\,\,
   \Phi_{j,\nu}=\sum_{i=j}^{d}\binom{i}{j}c_i\,p_\nu^{\,i-j}q_\nu^{\,d-i}=\frac{q_\nu^{d-j}}{j!}\,f^{(j)}\left(\frac{p_\nu}{q_\nu}\right) \in\mathbb{Z}.
   \end{equation}
  Put
  \begin{equation}\label{CC2}
 C_1 =
 \max_{0\le j\le  d}\,\,\max_{t\in[0,1]}\frac{|f^{(j)}(t)|}{j!} .
 \end{equation}
 Then
 \begin{equation}\label{c:Phisize}
|\Phi_{j,\nu}|\le C_1\,q_\nu^{d-j}\,\,\,\,(0\le j\le d).
\end{equation}
Also
\begin{equation}\label{Phimod}
\Phi_{0,\nu}=q_\nu^df\left(\frac{p_\nu}{q_\nu}\right)\equiv p_\nu^d,\,\,\,
\Phi_{1,\nu}  =q_\nu^{d-1}f'\left(\frac{p_\nu}{q_\nu}\right)
\equiv d\,p_\nu^{d-1},\,\,\,
\Phi_{2,\nu}
=\frac12q_\nu^{d-2}f''\left(\frac{p_\nu}{q_\nu}\right)
\equiv\binom d2p_\nu^{d-2}\pmod{q_\nu},
   \end{equation}
  as in  the sum  defining $ \Phi_{j,\nu}$ in (\ref{fvv}) all the terms with $i<d$ are divisible by $q_\nu$.
  In particular
   \begin{equation}\label{fvvv}
  ( \Phi_{0,\nu},q_\nu) = 1.
      \end{equation}
  
 \vskip+0.3cm
 Put $ P = d(d-1)(d-2) $ and take first partial quotients $a_1, a_2$ to satisfy
  \begin{equation}\label{begin}
  a_1\equiv 1 \pmod{P},\,\,\,\,\, a_2 \equiv 0 \pmod{P},\,\,\,\,\,
  f''(t) \neq 0 \,\,\,\text{for all}\,\,\, t \in \left(0, \frac{1}{a_1}\right]
  \end{equation}
  (we take into account that there exists $t_0$ such that $f''(t) \neq 0$ for all  $ t\in (0,t_0)$, as $f''$ has only finitely many zeros). 
   Then $\frac{p_1}{q_1},\frac{p_2}{q_2}$ and 
  $\Phi_{2,1}, \Phi_{2,2}$ are defined, and moreover
  $ \Phi_{2,1}\neq 0, \Phi_{2,2}\neq 0$.
   By (\ref{begin})   we see that
    \begin{equation}\label{begin11}
     |\Phi_{2,\nu}|\ge C_2\,q_\nu^{d-2},\,\,\,\,\,
     \text{where}
     \,\,\,\,\,
    C_2= \frac{1}{2}\cdot \min_{t \in \left[\frac{p_2}{q_2},\frac1{a_1}\right]}
  |f''(t)| >0.
  \end{equation}

   Take an integer $n$ with  $2^n > \max  (|  \Phi_{2,1}|, |\Phi_{2,2}|)$ and put
    \begin{equation}\label{K}
K=P^n = (d(d-1)(d-2))^n.
  \end{equation}
  
   \vskip+0.3cm
  Now we describe the inductive step of defining partial quotients.
  Assume that we have defined quotients $a_1,...,a_\nu, \nu \ge 2$ and 
  $$
  a_j \equiv 0\pmod{K}, 3\le j \le \nu.
  $$
  Define number 
      \begin{equation}\label{N}
      N_\nu=p_\nu^2\Phi_{2,\nu}-(d-1)p_\nu\Phi_{1,\nu}+(d-1)^2\Phi_{0,\nu}.
  \end{equation}
  By (\ref{Phimod}) we see that
  \begin{equation}\label{NN}
  N_\nu\equiv\binom{d-1}2p_\nu^d\pmod{q_\nu}.
  \end{equation}
  We define $D_\nu$ from the conditions
  \begin{equation}\label{D}
p_\nu\Phi_{0,\nu}D_\nu\equiv N_\nu\pmod{q_\nu^{d-1}},\,\,\,\,\,\,\,\,\,\,
-\frac{q_\nu^{d-1}}{2} <D_\nu\le\frac{q_\nu^{d-1}}{2} 
\end{equation}
(because of (\ref{fvvv}) by these conditions $D_\nu$ is defined uniquely). Now  for integer $a$
put
$$q_{\nu+1} = aq_\nu +q_{\nu-1}
$$
and define successively  three numbers
\begin{equation}\label{w}
\begin{array}{l}
\displaystyle w_{2,\nu}=\frac{K\left(D_\nu-(-1)^\nu(d-1)q_{\nu+1}\right)}{\Phi_{2,\nu}}
=
\frac{V_{2,\nu} a+U_{2,\nu}}{\Phi_{2,\nu}}
, \\[14pt]
\displaystyle w_{1,\nu}=\frac{(d-1)K-w_{2,\nu}\Phi_{1,\nu}-(-1)^\nu K    p_\nu q_{\nu+1}}{q_\nu^{d-1}}
=
\frac{V_{1,\nu} a+U_{1,\nu}}{\Phi_{2,\nu}\, q_\nu^{d-1}}, 
 \\[14pt]
\displaystyle w_{0,\nu}=\frac{K p_\nu-w_{2,\nu}\Phi_{0,\nu}-p_\nu q_\nu^{d-1}w_{1,\nu}}{q_\nu^{d}}
=
\frac{V_{0,\nu} a+U_{0,\nu}}{\Phi_{2,\nu}\, q_\nu^{d}},
\end{array}
\end{equation}
depending on $a$
(we took into account that  $\Phi_{2,\nu}\neq 0$
 because 
 $\frac{p_\nu}{q_\nu} \in \left[\frac{p_2}{q_2}, \frac{1}{a_1}\right]$ and 
 $ f''(t)$  has no zeros there by (\ref{begin})).
  Here integers $V_{j,\nu},U_{j,\nu}, j =0,1,2 $ do not  depend on $a$. For them we have explicit formulas
\begin{equation}\label{U}
\begin{array}{l}
\displaystyle 
 V_{2,\nu}  = (-1)^{\nu+1} K(d-1)q_{\nu} ,\,\,\,\,
U_{2,\nu}=K\,D_\nu+\frac{q_{\nu-1}}{q_\nu}V_{2,\nu},\\[10pt]
\displaystyle 
V_{1,\nu} 
=(-1)^{\nu}Kq_\nu ((d-1)\Phi_{1,\nu}-\Phi_{2,\nu}p_\nu),\,\,\,\,
U_{1,\nu}=K\left((d-1)\Phi_{2,\nu}-D_\nu\Phi_{1,\nu}\right)+\frac{q_{\nu-1}}{q_\nu}V_{1,\nu},\\[10pt]
\displaystyle 
V_{0,\nu} = (-1)^{\nu +1}K q_\nu ((d-1)(\Phi_{1,\nu}p_\nu -\Phi_{0,\nu}) -\Phi_{2,\nu}p_\nu^2),\,\,\,\,
U_{0,\nu}=K\left(D_\nu(\Phi_{1,\nu}p_\nu-\Phi_{0,\nu})-(d-2)p_\nu\Phi_{2,\nu}\right)+\frac{q_{\nu-1}}{q_\nu}V_{0,\nu}
\end{array}
\end{equation}
(each $V_{j,\nu}$ is divisible by $q_\nu$, so $U_{j,\nu}$ are integers). By direct calculation from (\ref{U})
and congruence (\ref{D})   for $D_\nu$  we deduce
\begin{equation}\label{det}
V_{1,\nu}U_{0,\nu}-V_{0,\nu}U_{1,\nu}=(-1)^\nu K^2q_\nu\,\Phi_{2,\nu}\left(p_\nu\Phi_{0,\nu}D_\nu-N_\nu\right)
\equiv0\pmod{q_\nu^d}
\end{equation}
(here $N_\nu$ is defined in (\ref{N})).

We are ready to formulate the conditions for choosing $a_{\nu+1}$. This partial quotient is chosen to be the minimal number from the set  $\{ a\in \frak{A}: a \ge q_\nu^{\tau-1}\}$, where
$$
\frak{A} =\{ a\in \mathbb{Z}:\,\,
  a \equiv 0\pmod{K},\,\,\,
\text{ all the numbers}\,\,\,  w_{0,\nu}, w_{1,\nu},  w_{2,\nu}\,\,\,
\text{are integers}\}.
$$

 \begin{lemma}\label{LA}   Let $\nu \ge 2$.
 Suppose that $a_1,a_2$ satisfy  {\rm (\ref{begin})} and that $a_3,...,a_\nu$ are divisible by $K$.
 Let 
 \begin{equation}\label{M}
 M_\nu = K\,|\Phi_{2,\nu}|\,q_\nu^{d-1} .
 \end{equation}
 Then for any integer $m$ there exists $a\in \frak{A}$ such that 
 $m \le a< m+M_\nu$.
  \end{lemma}

Lemma \ref{LA} shows that  $a_{\nu+1}$ satisfies conditions
  \begin{equation}\label{a0}
q_\nu^{\tau-1}\le a_{\nu+1}\le q_\nu^{\tau-1}+M_\nu,
\end{equation}
where $M_\nu$ is defined in (\ref{M}),
  \begin{equation}\label{a1}
a_{\nu+1} \equiv 0\pmod{K},
\end{equation}
where $K$ is defined in (\ref{K}),
and
  \begin{equation}\label{a2}
  w_{0,\nu}, w_{1,\nu},  w_{2,\nu}\,\,\,\text{are integers}.
\end{equation}

Now we prove Lemma \ref{LA}.
 \begin{proof}

We look for $a$ of the form $a=Kb$ with integer $b$.

\vskip+0.2cm

First of all we make two auxiliary observations.

\vskip+0.2cm
\noindent
{\bf 1.} We have $q_1=a_1\equiv1$, $q_2=a_1a_2+1\equiv1\pmod P$, and $q_{j+1}=a_{j+1}q_j+q_{j-1}\equiv q_{j-1}\pmod P$ when $a_{j+1}$ is
divisible by $K$. So $q_\nu\equiv1\pmod P$, and $q_\nu$ is coprime to $P$ and to $K$. By (\ref{Phimod}) and (\ref{U})
\begin{equation}\label{LAmodq}
\Phi_{2,\nu}\equiv\binom d2p_\nu^{d-2},\,\,\,\,\,\,\,\,\,\,\frac{V_{0,\nu}}{q_\nu}\equiv(-1)^{\nu+1}K\binom{d-1}2p_\nu^d\pmod{q_\nu} .
\end{equation}
As $\binom d2$ and $\binom{d-1}2$ divide $P$, the numbers $\Phi_{2,\nu}$ and $\frac{V_{0,\nu}}{q_\nu}$ are coprime to $q_\nu$.

\vskip+0.2cm
\noindent
{\bf 2.} Let $\ell$ be a prime divisor of $P$. As $a_3,...,a_\nu$ are divisible by $K$, we have 
$$
p_\nu \equiv p_i, \,\,\,\, q_\nu \equiv q_i \pmod{K} \,\,\,\,\,
\text{for some }\,\, i \in \{1,2\},
$$ and so $\Phi_{2,\nu}\equiv\Phi_{2,i}\pmod{\ell^n}$. As $\ell^n\ge 2^n>|\Phi_{2,i}|>0$, the number $\Phi_{2,i}$ is not divisible by $\ell^n$.
So $\Phi_{2,\nu}$ is not divisible by $\ell^n$.   As  prime $\ell$ can appear in the decomposition of $\Phi_{2,\nu}$ into primes only in a power less than  $n$, we conclude that 
 \begin{equation}\label{LAgcd}
\gcd\left(\Phi_{2,\nu},\,K^2(d-1)\right)\,\,\,\text{divides}\,\,\,K .
\end{equation}

Now we analyse when numbers $w_{j,\nu}$ are integers.

\vskip+0.2cm
\noindent
{\bf Number $w_{2,\nu}$.} By (\ref{w}) it is an integer if and only if
$$ 
V_{2,\nu}K\cdot b+U_{2,\nu}\equiv0\pmod{\Phi_{2,\nu}}.
$$
 This congruence has a
solution, because $\gcd(V_{2,\nu}K,\Phi_{2,\nu})=\gcd(K^2(d-1),\Phi_{2,\nu})$ divides $K$ by (\ref{LAgcd}), and $K$ divides $U_{2,\nu}$ by (\ref{U}). So
there exists $b_1$ such that $w_{2,\nu}\in\mathbb{Z}$ for all $b\equiv b_1\pmod{\Phi_{2,\nu}}$.

Let $w_{2,\nu}\in\mathbb{Z}$. Then by (\ref{w}) the numbers
$$
q_\nu^{d-1}w_{1,\nu}=\frac{V_{1,\nu}a+U_{1,\nu}}{\Phi_{2,\nu}}\,\,\,\,\,\,\,\text{and}\,\,\,\,\,\,\,q_\nu^dw_{0,\nu}=\frac{V_{0,\nu}a+U_{0,\nu}}{\Phi_{2,\nu}}
$$
are integers. As $\Phi_{2,\nu}$ is coprime to $q_\nu$, we conclude that 
\begin{equation}\label{LAcrit}
w_{1,\nu}\in\mathbb{Z}\,\Longleftrightarrow\,V_{1,\nu}a+U_{1,\nu}\equiv0\!\!\pmod{q_\nu^{d-1}},
\,\,\,\,\,\,\,\,\,\,\,\,
w_{0,\nu}\in\mathbb{Z}\,\Longleftrightarrow\,V_{0,\nu}a+U_{0,\nu}\equiv0\!\!\pmod{q_\nu^{d}}.
\end{equation}

\vskip+0.2cm
\noindent
{\bf Number $w_{0,\nu}$.} The number $U_{0,\nu}$ is divisible by $q_\nu$. Indeed, by (\ref{NN}), (\ref{D}) and (\ref{Phimod}) we have
$p_\nu D_\nu\equiv\binom{d-1}2\pmod{q_\nu}$,
and $p_\nu q_{\nu-1}\equiv(-1)^{\nu+1}\pmod{q_\nu}$ by (\ref{zero}). So by (\ref{U}), (\ref{Phimod}) and (\ref{LAmodq}) we get
$$
\frac{U_{0,\nu}}{K}\equiv D_\nu(d-1)p_\nu^d-(d-2)\binom d2p_\nu^{d-1}+(-1)^{\nu+1}\binom{d-1}2p_\nu^d\,q_{\nu-1}\equiv
 p_\nu^{d-1}\left(d\binom{d-1}2-(d-2)\binom d2\right)=0  \pmod{q_\nu}.
$$
So for $a=Kb$ the congruence $V_{0,\nu}a+U_{0,\nu}\equiv0\pmod{q_\nu^d}$ is equivalent to
$$
\frac{V_{0,\nu}}{q_\nu}K\cdot b+\frac{U_{0,\nu}}{q_\nu}\equiv0\pmod{q_\nu^{d-1}},
$$
where the coefficient $\frac{V_{0,\nu}}{q_\nu}K$ is coprime to $q_\nu$. So there exists $b_2$ such
that $w_{0,\nu}\in\mathbb{Z}$ for all $b\equiv b_2\pmod{q_\nu^{d-1}}$.

\vskip+0.2cm
\noindent
{\bf  Number $w_{1,\nu}$.} Let $w_{0,\nu}\in\mathbb{Z}$, that is $V_{0,\nu}a+U_{0,\nu}\equiv0\pmod{q_\nu^d}$. Then by (\ref{det}) the right hand side of
$$
V_{0,\nu}\left(V_{1,\nu}a+U_{1,\nu}\right)=V_{1,\nu}\left(V_{0,\nu}a+U_{0,\nu}\right)-\left(V_{1,\nu}U_{0,\nu}-V_{0,\nu}U_{1,\nu}\right)
$$
is divisible by $q_\nu^d$. As $V_{0,\nu}=q_\nu\cdot\frac{V_{0,\nu}}{q_\nu}$ with $\frac{V_{0,\nu}}{q_\nu}$ coprime to $q_\nu$, we get
$V_{1,\nu}a+U_{1,\nu}\equiv0\pmod{q_\nu^{d-1}}$, that is
$w_{1,\nu}\in\mathbb{Z}$ by (\ref{LAcrit}). We see that $w_{0,\nu}\in \mathbb{Z}$  ensures $w_{1,\nu}\in \mathbb{Z}$  automatically.

\vskip+0.2cm
 Finally,  as  $\Phi_{2,\nu}$ and $q_\nu^{d-1}$ are coprime, by the Chinese Remainder Theorem there exists $b_0$ such that every
$$
b\equiv b_0\pmod{|\Phi_{2,\nu}|\,q_\nu^{d-1}}
$$
satisfies $b\equiv b_1\pmod{\Phi_{2,\nu}}$ and $b\equiv b_2\pmod{q_\nu^{d-1}}$. For such $b$ the number $a=Kb$ belongs
to $\frak A$. These numbers $a$ form an arithmetic progression with difference $K|\Phi_{2,\nu}|q_\nu^{d-1}=M_\nu$, and so every
interval $m\le a<m+M_\nu$ contains one of them.
 \end{proof}

  The next statement is a detailed version of Theorem \ref{T2}.

\begin{thm}\label{T4}
Let $\tau >2d-2$.
Let $\xi$ be defined by continued fraction {\rm (\ref{o})} with partial quotients satisfying
{\rm(\ref{begin})} and {\rm(\ref{a0}), (\ref{a1}), (\ref{a2})} for all $\nu\geqslant2$. Then $\xi$ is transcendental and
 \begin{equation}\label{tt4}
\hat\omega^*[f]\geqslant\omega_d(\tau)=\min\left\{\frac{\tau^2}{(d-1)(\tau+1)},\,\,\,\frac{\tau+1}{d},\,\,\,2+\frac{d-2}{\tau-d+2},\,\,\,
2+\frac{2}{\tau}-\frac{d-1}{\tau^2}\right\}.
\end{equation}
\end{thm}

 Theorem \ref{T2} follows from Theorem \ref{T4} immediately.
 
 \vskip+0.2cm
 Let $\xi$ satisfy the conditions of Theorem \ref{T4}. Then for all $\nu$ large enough
 by formulas (\ref{zero}),  (\ref{c:Phisize}), (\ref{M}) and  (\ref{a0})
  one has
\begin{equation}\label{growth}
q_\nu^{\tau-1}\leqslant a_{\nu+1}\leqslant2q_\nu^{\tau-1},\,\,\,\,\,\,\,\,\,\,
q_\nu^{\tau}\leqslant q_{\nu+1}\leqslant3q_\nu^{\tau},\,\,\,\,\,\,\,\,\,\,q_{\nu-1}\leqslant q_\nu^{1/\tau}.
\end{equation}
In particular $\xi$ is transcendental by Roth's theorem.

 \vskip+0.2cm
 In the rest of the paper we give a proof of inequality (\ref{tt4}) of Theorem \ref{T4}.

\subsection{Upper bound for  $D_\nu$}

Number $D_\nu$ is defined in (\ref{D}) as a residue modulo $q_\nu^{d-1}$, so a priori we know only that $|D_\nu|\le \frac{q_\nu^{d-1}}{2}$.
In fact it is much smaller. The reason is the following. For an integer $x$ coprime to $q_\nu$ we denote by $x^*$ an integer with 
$$
x\,x^*\equiv 1\pmod{q_\nu^{d-1}}.
$$
By (\ref{D}) we have 
$$
D_\nu \equiv N_\nu\, p_\nu^*\,\Phi_{0,\nu}^* \pmod{q_\nu^{d-1}},
$$
 and so we need the inverse $p_\nu^*$ of $p_\nu$ modulo $q_\nu^{d-1}$.
From the equality $q_\nu p_{\nu-1}-p_\nu q_{\nu-1} = (-1)^\nu$ we see that $p_\nu\cdot (-1)^{\nu+1}q_{\nu-1} = 1- u q_\nu$ with $u = (-1)^\nu p_{\nu-1}$, and therefore we can take
$$
p_\nu^* = (-1)^{\nu+1}q_{\nu-1} \left( 1 + uq_\nu+ (uq_\nu)^2+...+(uq_\nu)^{d-2}\right) .
$$
So $p_\nu^*$ is a number of size $\asymp q_{\nu-1}(q_\nu q_{\nu-1})^{d-2}$, which is much smaller than the modulus $q_\nu^{d-1}$ as soon as $q_{\nu-1}^{d-1}$ is
small compared with $q_\nu$. The same happens with the whole product $N_\nu\, p_\nu^*\,\Phi_{0,\nu}^*$, as modulo $q_\nu^{d-1}$ it can be written as a sum $\sum_{m=0}^{d-2}s_mq_\nu^m$ with coefficients $|s_m|\ll q_{\nu-1}^{m+1}$.
Thus $D_\nu$ is small because of the general identity for convergents together with the fast growth of the denominators $q_\nu$.

 We begin with an auxiliary statement about lifting modulo $q^m$.

\begin{lemma}\label{e:Laux}
Let $q, r\ge 1$ and $n\ge 0$ be integers and $\gamma \ge 1$. Suppose that integers $g_m, h_m$, $0\le m\le n$ satisfy
$$
g_0 = 1,\,\,\,\,\,\,\, |g_m|\le \gamma r^m,\,\,\,\,\,\,\, |h_m|\le \gamma r^m, \,\,\,\,\,\,\, 0\le m \le n.
$$
Then there exist integers $s_0,s_1,...,s_n$ such that
\begin{equation}\label{e:aux}
|s_m|\le \gamma(1+\gamma)^m r^m
\,\,\,\,\,\,\,\text{and}\,\,\,\,\,\,\,
\left(\sum_{m=0}^n g_mq^m\right)\cdot \left(\sum_{m=0}^n s_mq^m\right)\equiv \sum_{m=0}^n h_mq^m\pmod{q^{n+1}}.
\end{equation}
\end{lemma}
\begin{proof}
Define $s_0,s_1,...,s_n$ successively by
$$
s_0 = h_0,\,\,\,\,\,\,\,\,\,
s_m = h_m - \sum_{k=1}^m g_k s_{m-k},\,\,\,\, 1\le m \le n.
$$
Then for every $m\le n$ the coefficient of $q^m$ in the product from (\ref{e:aux}) is equal to $\sum_{k=0}^m g_ks_{m-k} = h_m$, and all
other terms of the product are divisible by $q^{n+1}$. This gives the congruence. The inequality is proven by induction as
$$
|s_m|\le \gamma r^m + \sum_{k=1}^m \gamma r^k\cdot \gamma (1+\gamma)^{m-k}r^{m-k} =\gamma r^m\left( 1+\gamma \sum_{j=0}^{m-1}(1+\gamma)^j\right) = \gamma (1+\gamma)^mr^m.
$$
\end{proof}

\begin{lemma}\label{e:LD}
Let $\xi$ satisfy the conditions of Theorem \ref{T4}. 
Then there exists a constant $C_3$ depending only on $f$ such that for all $\nu$ large enough
$$
|D_\nu|\le C_3\, q_\nu^{d-2}q_{\nu-1}^{d-1}.
$$
\end{lemma}

\begin{proof}
Put
$$
e_\nu = (-1)^{\nu+1}q_{\nu-1},\,\,\,\,\,\,\, u_\nu = (-1)^\nu p_{\nu-1},\,\,\,\,\,\,\,\text{so that}\,\,\,\,\,\,\, |e_\nu|\le q_{\nu-1},\,\,\,\, |u_\nu|\le q_{\nu-1}.
$$
The equality $q_\nu p_{\nu-1}-p_\nu q_{\nu-1} = (-1)^\nu$ from (\ref{zero}) can be written as
\begin{equation}\label{e:pe}
p_\nu e_\nu = 1 - u_\nu q_\nu .
\end{equation}
In particular $e_\nu$ is coprime to $q_\nu$.

By the definition (\ref{fvv}) of $\Phi_{j,\nu}$ and by (\ref{e:pe}) for $j = 0,1,2$ we have
\begin{equation}\label{e:ePhi}
e_\nu^{d-j}\Phi_{j,\nu} = \sum_{i=j}^d \binom{i}{j}c_i\, (p_\nu e_\nu)^{i-j}(q_\nu e_\nu)^{d-i}
=
\sum_{i=j}^d \binom{i}{j}c_i\, (1-u_\nu q_\nu)^{i-j}(e_\nu q_\nu)^{d-i} .
\end{equation}
Now we represent 
  the right hand side of (\ref{e:ePhi}) as  a sum
  $$
  \sum_{m} h_m q_\nu^m,
  $$ 
  where $h_m$ is a sum of at most $\gamma$ products of $m$ factors, each factor being equal to $\pm u_\nu$ or to $\pm e_\nu$,
with a certain $\gamma$ depending only on $f$. So  we get 
$$
|h_m|\le \gamma q_{\nu-1}^m.
$$
 By (\ref{e:pe}) and (\ref{N}) the same is true for the numbers
$$
G_\nu =e_\nu^{d+1}\, p_\nu \Phi_{0,\nu}= (1-u_\nu q_\nu)\cdot e_\nu^d\Phi_{0,\nu} 
\,\,\,\,\,
\text{and}
\,\,\,\,\,
H_\nu = e_\nu^d N_\nu =
(1-u_\nu q_\nu)^2\cdot e_\nu^{d-2}\Phi_{2,\nu} - (d-1)(1-u_\nu q_\nu)\cdot e_\nu^{d-1}\Phi_{1,\nu}+(d-1)^2e_\nu^d\Phi_{0,\nu}.
$$
These numbers can  be written as 
$$
G_\nu
= \sum_m g_m q_\nu^m,\,\,\,\,\,\,\,\,
H_\nu  
=\sum_m h_m q_\nu^m,\,\,\,\,\,
\text{
with}
\,\,\,\,\,\,\,\,
 |g_m|\le \gamma q_{\nu-1}^m,\,\,\, |h_m|\le \gamma q_{\nu-1}^m.
 $$
  Here  $\gamma$ depends only on $f$, and  $g_0 = c_d =1$.

Now we apply Lemma \ref{e:Laux} with $q = q_\nu$, $r= q_{\nu-1}$ and $n = d-2$ (terms with $q_\nu^m, m\ge d-1$ in $G_\nu$ and $H_\nu$ are divisible by $q_\nu^{d-1}$ and
may be omitted). We get an integer
$$
S_\nu = \sum_{m=0}^{d-2}s_m q_\nu^m
\,\,\,\,\,\,\,\text{with}\,\,\,\,\,\,\,
G_\nu S_\nu \equiv H_\nu \pmod{q_\nu^{d-1}}
$$
and
$$
|S_\nu| \le \gamma(1+\gamma)^{d-2}\sum_{m=0}^{d-2}(q_\nu q_{\nu-1})^m\le 2\gamma(1+\gamma)^{d-2}(q_\nu q_{\nu-1})^{d-2}.
$$
The congruence here means that $e_\nu^{d+1}\, p_\nu \Phi_{0,\nu} S_\nu\equiv e_\nu^d N_\nu \pmod{q_\nu^{d-1}}$. As $e_\nu$ is coprime to $q_\nu$, for the integer
$$
D_\nu' = e_\nu S_\nu
$$
we get
\begin{equation}\label{e:Dprime}
p_\nu\Phi_{0,\nu}D_\nu' \equiv N_\nu \pmod{q_\nu^{d-1}},\,\,\,\,\,\,\,\,\,\,
|D_\nu'|\le C_3\, q_\nu^{d-2}q_{\nu-1}^{d-1},\,\,\,\,\,\,\, C_3 = 2\gamma(1+\gamma)^{d-2}.
\end{equation}
By (\ref{a0}) we have $q_\nu \ge q_{\nu-1}^\tau$ with $\tau > d-1$. So $C_3q_{\nu-1}^{d-1}<\frac{q_\nu}{2}$ for all $\nu$ large enough, and then
$|D_\nu'|<\frac{q_\nu^{d-1}}{2}$. We see that $D_\nu'$ satisfies both conditions from (\ref{D}). As these conditions define $D_\nu$ uniquely, we have $D_\nu = D_\nu'$,
and the lemma follows from (\ref{e:Dprime}).
\end{proof}

\subsection{Vectors $A_\nu,B_\nu$ and $F_\nu$}

 As in Section \ref{s2},
 for vector $X=(x_0,x_1,x_2)\in\mathbb{R}^3$ we denote
$$
|X|=\max_{j=0,1,2}|x_j|,\,\,\,\,\,\,\,L(X)=|x_0+x_1\xi+x_2f(\xi)| 
$$
 and
 put $A_\nu=(-p_\nu,q_\nu,0)$,
but now $B_\nu$  will denote the vector
\begin{equation}\label{c:AB}
B_\nu=\left(p_\nu\Phi_{1,\nu}-\Phi_{0,\nu},\,\,-q_\nu\Phi_{1,\nu},\,\,q_\nu^d\right) \in \mathbb{Z}^3.
\end{equation}
For $f(t)=t^3$ this is the vector $B_\nu=(2p_\nu^3,-3q_\nu p_\nu^2,q_\nu^3)$ from (\ref{AB}).

Analogously to Lemma \ref{L1} for $A_\nu$ we have
\begin{equation}\label{XA}
|A_\nu|=q_\nu,\,\,\,\,\,\,\,\,L(A_\nu)<\frac1{q_{\nu+1}},\,\,\,\,\nu \ge 1.
\end{equation}
For  the new vector $B_\nu$ now we formulate an analog of Lemma \ref{L2}.

\begin{lemma}\label{XL2}
For every $\nu\ge1$ one has $B_\nu\in\mathbb{Z}^3$ and 
$$
|B_\nu|\le 2C_1q_\nu^d,\,\,\,\,\,L(B_\nu)\le 2C_1\frac{q_\nu^{d-2}}{q_{\nu+1}^2},
$$
where 
$C_1$ is defined in {\rm(\ref{CC2})}.
 
\end{lemma}

\begin{proof}
The bound for $|B_\nu|$ follows from (\ref{c:Phisize}). By    (\ref{c:AB}), (\ref{fvv}), (\ref{c:Phisize}) and (\ref{11})
we get the  upper bound
$$
L(B_\nu) =
|
p_\nu\Phi_{1,\nu}-\Phi_{0,\nu} -q_\nu\Phi_{1,\nu}\xi +q_\nu^d f(\xi)
|
\le
\sum_{j=2}^d|\Phi_{j,\nu}|\cdot|q_\nu\xi-p_\nu|^j\le
\sum_{j=2}^d\frac{C_1q_\nu^{d-j}}{q_{\nu+1}^j}
=\frac{C_1q_\nu^{d-2}}{q_{\nu+1}^2}\sum_{j=2}^d\left(\frac{1}{q_\nu q_{\nu+1}}\right)^{j-2}\le
\frac{2C_1q_\nu^{d-2}}{q_{\nu+1}^2}.
$$
\end{proof}

Vectors $A_\nu$ and $B_\nu$ are not enough for our proof, because  between $|A_\nu|=q_\nu$ and $|B_\nu|\asymp q_\nu^d$ the gap is too large. The
next vector is a linear combination of $A_\nu$ and $B_\nu$ with rational coefficients, which has integer
coordinates and smaller norm than $B_\nu$. Put
\begin{equation}\label{c:F}
F_\nu=\frac{\kappa_\nu}{q_\nu}A_\nu+\frac{q_{\nu-1}^{d-1}}{q_\nu}B_\nu,\,\,\,\,\,\,\,\,\,\,\text{where}\,\,\,\,\,\kappa_\nu=(-1)^{(\nu+1)(d-1)}(d-1).
\end{equation}

\begin{lemma}\label{c:LF}
For every $\nu\ge2$ vector $F_\nu$ has integer coordinates
\begin{equation}\label{c:Fcoord}
F_\nu=\left(\frac{q_{\nu-1}^{d-1}\left(p_\nu\Phi_{1,\nu}-\Phi_{0,\nu}\right)-\kappa_\nu p_\nu}{q_\nu},\,\,\,\,\kappa_\nu-q_{\nu-1}^{d-1}\Phi_{1,\nu},\,\,\,\,(q_\nu q_{\nu-1})^{d-1}\right),
\end{equation}
and
$$
|F_\nu|\le 3C_1(q_\nu q_{\nu-1})^{d-1},\,\,\,\,\,\,\,\,\,\,
L(F_\nu)\le\frac{1}{q_\nu q_{\nu+1}}\left(d-1+2C_1\,\frac{q_{\nu-1}^{d-1}q_\nu^{d-2}}{q_{\nu+1}}\right),
$$
where $C_1$ is defined in {\rm(\ref{CC2})}.
\end{lemma}

\begin{proof}
Formula (\ref{c:Fcoord}) follows from the definition of $A_\nu$ and  (\ref{c:AB}) and (\ref{c:F}). The second and the third coordinates here are
integers. For the first coordinate we use (\ref{Phimod}) and then the last formula from the first row  from (\ref{zero}).
In such a way we get
$$
\begin{array}{l}
q_{\nu-1}^{d-1}\left(p_\nu\Phi_{1,\nu}-\Phi_{0,\nu}\right)-\kappa_\nu p_\nu\equiv
q_{\nu-1}^{d-1}\left(d\,p_\nu^d-p_\nu^d\right)-\kappa_\nu p_\nu=\\[4pt]
\,\,\,\,\,\,\,\,\,\,=(d-1)p_\nu\left((p_\nu q_{\nu-1})^{d-1}-(-1)^{(\nu+1)(d-1)}\right)\equiv0 \pmod{q_\nu} .
\end{array}
$$
So $F_\nu\in\mathbb{Z}^3$. The bound for $|F_\nu|$ follows from (\ref{c:Fcoord}) and (\ref{c:Phisize}). 
For the linear form $L(F_\nu)$ from the definitions, the bound  (\ref{XA}) for $L(A_\nu)$ and Lemma \ref{XL2}
for $B_\nu$ we get 
$$
L(F_\nu) = \frac{1}{q_\nu}\left| \kappa_\nu (q_\nu\xi-p_\nu) + q_{\nu-1}^{d-1}\sum_{j=2}^d\Phi_{j,\nu}\,(q_\nu\xi-p_\nu)^j\right|
\le \frac{(d-1)L(A_\nu)+q_{\nu-1}^{d-1}L(B_\nu)}{q_\nu}<
$$
$$
<\frac{d-1}{q_\nu q_{\nu+1}} + \frac{q_{\nu-1}^{d-1}}{q_\nu}\cdot\frac{2C_1q_\nu^{d-2}}{q_{\nu+1}^2}
=\frac{1}{q_\nu q_{\nu+1}}\left( d-1 +\frac{2C_1\,q_{\nu-1}^{d-1}q_\nu^{d-2}}{q_{\nu+1}}\right)
 .
$$
Lemma is proven.
\end{proof}

\subsection{Vector $W_\nu$}

Let $\xi$ satisfy the conditions of Theorem \ref{T4}. For $\nu\geqslant2$ we consider the vector
$$
W_\nu=(w_{0,\nu},w_{1,\nu},w_{2,\nu})
$$
whose coordinates are the numbers (\ref{w}) with $a=a_{\nu+1}$. By (\ref{a2}) we have $W_\nu\in\mathbb{Z}^3$. The meaning of the
formulas (\ref{w}) is the following. By (\ref{fvv}) for any vector $X=(x_0,x_1,x_2)$ we have
$$
q_\nu^d L(X) =
\left|
\left(x_0q_\nu^d+x_1p_\nu q_\nu^{d-1}+x_2\Phi_{0,\nu}\right)+\left(x_1q_\nu^{d-1}+x_2\Phi_{1,\nu}\right)(q_\nu \xi-p_\nu)+x_2\sum_{j=2}^d\Phi_{j,\nu}\,(q_\nu \xi-p_\nu)^j \right|.
$$
The second and the third formulas from (\ref{w}) mean that for $X=W_\nu$ the first two brackets are equal to
$Kp_\nu$ and to $K\left(d-1-(-1)^\nu p_\nu q_{\nu+1}\right)$ and so
\begin{equation}\label{c:WTaylor}
q_\nu^d L(W_\nu) 
=
\left|
Kp_\nu+K\left(d-1-(-1)^\nu p_\nu q_{\nu+1}\right)(q_\nu \xi-p_\nu)+w_{2,\nu}\sum_{j=2}^d\Phi_{j,\nu}\,(q_\nu \xi-p_\nu)^j 
\right|
.
\end{equation}

\begin{lemma}\label{LW}
For all $\nu$ large enough one has
\begin{equation}\label{c:LWb}
w_{2,\nu}\neq0,\,\,\,\,\,\,\,\,\,\,
|W_\nu|\le C_4\frac{q_{\nu+1}}{q_\nu^{d-2}},\,\,\,\,\,\,\,
L(W_\nu)\le C_5\frac{q_{\nu-1}^{d-1}}{q_\nu^2q_{\nu+1}^2},
\end{equation}
 {where}
 $$
C_4 = K(d+1)+\frac{2C_1Kd}{C_2},\,\,\,\,
C_5 = KC_3+\frac{2C_1Kd}{C_2}+1.
$$
\end{lemma}

\begin{proof}
By Lemma \ref{e:LD} and (\ref{growth}) we have $|D_\nu|\le C_3q_\nu^{d-2}q_{\nu-1}^{d-1}\le q_\nu^{d-1}<q_{\nu+1}$ for $\nu$ large enough. So from the
first formula from (\ref{w}) and the lower bound (\ref{begin11}) we get
$$
w_{2,\nu}\neq0,\,\,\,\,\,\,\,\,\,\,|w_{2,\nu}|\le \frac{K\left( |D_\nu| + (d-1)q_{\nu+1}\right)}{|\Phi_{2,\nu}|}\le\frac{K\,d\,q_{\nu+1}}{|\Phi_{2,\nu}|}\le \frac{Kd}{C_2}\cdot\frac{q_{\nu+1}}{q_\nu^{d-2}}.
$$
Then by the second and the third formulas from (\ref{w}) and by (\ref{c:Phisize}) we get upper bounds
$$
|w_{1,\nu}|\le\frac{K\,d\,q_\nu q_{\nu+1}+|w_{2,\nu}|\cdot C_1q_\nu^{d-1}}{q_\nu^{d-1}}\le 
\left(Kd+\frac{C_1Kd}{C_2}\right)\cdot 
 \frac{q_{\nu+1}}{q_\nu^{d-2}},
$$
$$
|w_{0,\nu}|\le\frac{Kq_\nu+|w_{2,\nu}|\cdot C_1q_\nu^d+|w_{1,\nu}|\,q_\nu^d}{q_\nu^d}\le C_4\frac{q_{\nu+1}}{q_\nu^{d-2}}.
$$
This gives the bound for $|W_\nu|$ as $C_4$ is the maximal value among the constants in the right hand sides.

\vskip+0.2cm
Now we get the upper bound for $
L(W_\nu)$.
Consider polynomials
$$
A(T) = q_\nu T-p_\nu,\,\,\,\,\,\,\,
A_+(T) = q_{\nu+1} T-p_{\nu+1},
$$
$$
B(T) = \sum_{j=2}^d\Phi_{j,\nu}\,(q_\nu T-p_\nu)^j
\,\,\,\,\,
\text{and}
\,\,\,\,\,
W(T) = w_{0,\nu}+w_{1,\nu}T+w_{2,\nu}f(T).
$$
By (\ref{fvv}) and by the second and the third formulas from (\ref{w}) we have
$$
q_\nu^d W(T) = \left(w_{0,\nu}q_\nu^d+w_{1,\nu}p_\nu q_\nu^{d-1}+w_{2,\nu}\Phi_{0,\nu}\right) + \left(w_{1,\nu}q_\nu^{d-1}+w_{2,\nu}\Phi_{1,\nu}\right)A(T)+w_{2,\nu}B(T)=
$$
$$
=Kp_\nu + K\left(d-1-(-1)^\nu p_\nu q_{\nu+1}\right)A(T)+w_{2,\nu}B(T).
$$
As $q_{\nu+1}A(T) - (-1)^\nu = q_\nu A_+(T)$ by (\ref{zero}), this means that
\begin{equation}\label{v:xyz}
W(T) = xA(T)+yB(T)+zA_+(T)
\,\,\,\,\,\,\,
\text{with}
\,\,\,\,\,\,\,
x = \frac{K(d-1)}{q_\nu^d},\,\,\,\,\,
y = \frac{w_{2,\nu}}{q_\nu^d},\,\,\,\,\,
z = (-1)^{\nu+1}\frac{Kp_\nu}{q_\nu^{d-1}}
\end{equation}
(this is an analog of formulas (\ref{coe}) from Section 2). Notice that by (\ref{zero}) we have equalities
$$
A\left(\frac{p_{\nu+1}}{q_{\nu+1}}\right) =\frac{(-1)^\nu }{q_{\nu+1}},\,\,\,\,\,\,\,
A_+\left(\frac{p_{\nu+1}}{q_{\nu+1}}\right) =0
\,\,\,\,\,\,\,
  \text{and}
  \,\,\,\,\,\,\,
B\left(\frac{p_{\nu+1}}{q_{\nu+1}}\right) =\frac{\Phi_{2,\nu}}{q_{\nu+1}^2}+\sum_{j=3}^d\frac{(-1)^{\nu j}\Phi_{j,\nu}}{q_{\nu+1}^j}.
$$
So by (\ref{v:xyz}) we obtain
$$
W\left(\frac{p_{\nu+1}}{q_{\nu+1}}\right) = xA\left(\frac{p_{\nu+1}}{q_{\nu+1}}\right)
  + y B\left(\frac{p_{\nu+1}}{q_{\nu+1}}\right)
  =
  \frac{(-1)^\nu K(d-1)q_{\nu+1} + w_{2,\nu}\Phi_{2,\nu}}{q_\nu^dq_{\nu+1}^2} +
  \frac{w_{2,\nu}}{q_\nu^d}\sum_{j=3}^d\frac{(-1)^{\nu j}\Phi_{j,\nu}}{q_{\nu+1}^j}.
$$
By the first formula from (\ref{w}) we have $w_{2,\nu}\Phi_{2,\nu}+(-1)^\nu K(d-1)q_{\nu+1}=KD_\nu$, and after this cancellation we get
\begin{equation}\label{v:value}
W\left(\frac{p_{\nu+1}}{q_{\nu+1}}\right)=\frac{KD_\nu}{q_\nu^dq_{\nu+1}^2}+\frac{w_{2,\nu}}{q_\nu^d}\sum_{j=3}^d\frac{(-1)^{\nu j}\Phi_{j,\nu}}{q_{\nu+1}^j}.
\end{equation}
By Lemma \ref{e:LD} for the first summand here we have
$$
\left|\frac{KD_\nu}{q_\nu^dq_{\nu+1}^2}\right|\le KC_3\frac{q_{\nu-1}^{d-1}}{q_\nu^2q_{\nu+1}^2}.
$$
For the second summand we use (\ref{c:Phisize}) and the bound $|w_{2,\nu}|\le \frac{Kd\,q_{\nu+1}}{|\Phi_{2,\nu}|}\le \frac{Kd}{C_2}\cdot\frac{q_{\nu+1}}{q_\nu^{d-2}}$ obtained above. In such a way we obtain the estimate
$$
\left|\frac{w_{2,\nu}}{q_\nu^d}\sum_{j=3}^d\frac{(-1)^{\nu j}\Phi_{j,\nu}}{q_{\nu+1}^j}\right|\le
\frac{Kd}{C_2}\cdot\frac{q_{\nu+1}}{q_\nu^{d-2}}\cdot\frac1{q_\nu^d}\cdot\frac{2C_1q_\nu^{d-3}}{q_{\nu+1}^3}\le\frac{2C_1Kd}{C_2}\cdot\frac{1}{q_\nu^2q_{\nu+1}^2}.
$$
So for all $\nu$ large enough we get
\begin{equation}\label{v:P2}
 \left| W\left(\frac{p_{\nu+1}}{q_{\nu+1}}\right)\right| \le \left(KC_3+\frac{2C_1Kd}{C_2}\right)\frac{q_{\nu-1}^{d-1}}{q_\nu^2q_{\nu+1}^2}.
\end{equation}
On the other hand, for all $\nu$ large enough we have
\begin{equation}\label{v:P1}
 \left| W(\xi) - W\left(\frac{p_{\nu+1}}{q_{\nu+1}}\right)\right|\le
 C_6 |W_\nu| \cdot
  \left|\xi - \frac{p_{\nu+1}}{q_{\nu+1}}\right|\le
  \frac{C_4C_6}{q_\nu^{d-2} q_{\nu+2}}
 \le
  \frac{1}{q_\nu^2q_{\nu+1}^2},
  \,\,\,\,\,\,\,\,
  C_6 = 1+\max_{0\le t\le 1}|f'(t)|
\end{equation}
(we used (\ref{11}), the upper bound for $|W_\nu|$ and the inequalities $q_{\nu+2}\ge q_{\nu+1}^\tau$, $\tau>4$ from (\ref{growth})).
Now combining (\ref{v:P2}) and (\ref{v:P1}) we get the inequality for $L(W_\nu) = |W(\xi)|$.
\end{proof}

\subsection{Proof of Theorem \ref{T4}}

 We should prove (\ref{tt4}). If $\omega_d(\tau)\le 2$, there is nothing to prove. Let
$2<\omega<\omega_d(\tau)$. We prove that
\begin{equation}\label{c:goal}
\lim_{t\to\infty}t^{\omega}\psi^*(t)=0.
\end{equation}
Then $\hat\omega^*[f]\geqslant\omega$, and this gives Theorem \ref{T4}. Let $C_7 = \max (3C_1, C_4)\ge1$ be the maximum of the constants from the
bounds for $|F_\nu|$, $|B_\nu|$, $|W_\nu|$ in Lemmas \ref{XL2}, \ref{c:LF}, \ref{LW}.
Then
$$
|F_\nu| \le C_7 (q_\nu q_{\nu-1})^{d-1}, \,\,\,\, |B_\nu| \le C_7q_\nu^d, \,\,\,\, |W_\nu | \le C_7 
\frac{
q_{\nu+1}}{q_\nu^{d-2}}.
$$
 By (\ref{growth}) for all $\nu$ large enough
we have
$$
q_\nu<C_7(q_\nu q_{\nu-1})^{d-1}<C_7q_\nu^d<C_7\frac{q_{\nu+1}}{q_\nu^{d-2}}<q_{\nu+1}.
$$
Let $q_\nu\le t<q_{\nu+1}$. In each of the vectors $A_\nu,F_\nu,B_\nu,W_\nu$ the second or the third coordinate is
not equal to zero. 
We  consider four cases and use the inequalities (\ref{growth}),
(\ref{XA}) and all the inequalities  for linear forms from  Lemmas \ref{XL2}, \ref{c:LF}, \ref{LW}.

\medskip
\noindent
{\bf 1.} For $q_\nu\le t<C_7(q_\nu q_{\nu-1})^{d-1}$ we have $\psi^*(t)\le L(A_\nu)$, and from  (\ref{XA}) we get
$$
t^\omega\psi^*(t)\le \left(C_7(q_\nu q_{\nu-1})^{d-1}\right)^\omega\frac{1}{q_{\nu+1}}\le C_8\,q_\nu^{(d-1)\left(1+\frac{1}{\tau}\right)\omega-\tau},\,\,\,\,\,\text{with}\,\,\,\,\, C_8 = C_7^\omega.
$$
The exponent is negative, because $\omega<\dfrac{\tau^2}{(d-1)(\tau+1)}$.

\medskip
\noindent
{\bf 2.} For $C_7(q_\nu q_{\nu-1})^{d-1}\le t<C_7q_\nu^d$ we have $\psi^*(t)\le L(F_\nu)$. From Lemma  \ref{c:LF} we 
deduce upper bound
$
L(F_\nu)\le \frac{d}{q_\nu q_{\nu+1}}$ for large $\nu$, and so 
$$
t^\omega\psi^*(t)\le \left(C_7q_\nu^d\right)^\omega\frac{d}{q_\nu q_{\nu+1}}\le C_9\,q_\nu^{d\omega-1-\tau}
,\,\,\,\,\,\text{with}\,\,\,\,\, C_9 = d\,C_7^\omega.
$$
The exponent is negative, because $\omega<\dfrac{\tau+1}{d}$.

\medskip
\noindent
{\bf 3.} For $C_7q_\nu^d\leq t<C_7\dfrac{q_{\nu+1}}{q_\nu^{d-2}}$ we have $\psi^*(t)\le L(B_\nu)$, and from Lemma \ref{XL2} we
get
$$
t^\omega\psi^*(t)\le  \left(C_7\frac{q_{\nu+1}}{q_\nu^{d-2}}\right)^\omega\cdot 2C_1\frac{q_\nu^{d-2}}{q_{\nu+1}^2}
=2C_1C_7^\omega\frac{q_{\nu+1}^{\omega-2}}{q_\nu^{(d-2)(\omega-1)}}\le C_{10}\,q_\nu^{\tau(\omega-2)-(d-2)(\omega-1)}
,\,\,\,\,\,\text{with}\,\,\,\,\, C_{10} = 2\cdot 3^{\omega-2}C_1C_7^\omega.
$$
The exponent is negative, because $\omega(\tau-d+2)<2\tau-d+2$, that is $\omega<2+\dfrac{d-2}{\tau-d+2}$.

\medskip
\noindent
{\bf 4.} For $C_7\dfrac{q_{\nu+1}}{q_\nu^{d-2}}\le t <q_{\nu+1}$ we have $\psi^*(t)\le L(W_\nu)$, and from Lemma \ref{LW} we get
$$
t^\omega\psi^*(t)\le q_{\nu+1}^\omega\cdot C_5\frac{q_{\nu-1}^{d-1}}{q_\nu^2q_{\nu+1}^2}
=C_5\frac{q_{\nu+1}^{\omega-2}q_{\nu-1}^{d-1}}{q_\nu^2}\le C_{11}\,q_\nu^{\tau(\omega-2)+\frac{d-1}{\tau}-2},
\,\,\,\,\,\text{with}\,\,\,\,\, C_{11} = 3^{\omega-2}C_5.
$$
The exponent is negative, because $\omega<2+\dfrac{2}{\tau}-\dfrac{d-1}{\tau^2}$.

\medskip
The four inequalities for $\omega$ are just the inequalities $\omega<\omega_d(\tau)$ for the four numbers from
the definition of $\omega_d(\tau)$. So in all the cases $t^\omega\psi^*(t)\le C_{12} q_\nu^{-\delta},
C_{12} = \max(C_8,C_9,C_{10},C_{11})$
 with some $\delta>0$, and (\ref{tt4}) follows. $\Box$

\vskip+0.3cm
 {\bf AI declaration.}
This paper was prepared with the assistance of Claude, an artificial intelligence system
developed by Anthropic.
The author is grateful for  this assistance.
In particular, the main idea of the construction, namely to link the successive blocks of
approximations generated by the intrinsic rational points
$\left(\frac{p_\nu}{q_\nu}, f\left(\frac{p_\nu}{q_\nu}\right)\right)$ of the curve
by means of the vectors $W_\nu$, was discovered by the AI system in the course of a dialogue with the author.
Claude also assisted in the calculations and 
in checking the proofs. The author takes full responsibility for the  correctness of the paper.

\vskip+0.3cm
 {\bf Acknowledgements.}
 This research is supported by 
  by Austrian Science Fund (FWF), Forschungsprojekt PAT1961524.

\end{document}